
PAPER:
%
%
%
\input amstex
\documentstyle{amsppt}
\magnification=1200
%
%
%
\hoffset=0.4cm\voffset=0.4cm
\hsize=13.5cm\vsize=18.5cm
\baselineskip=13pt plus 0.2pt
\parskip=4pt
%
%
%
\def\R{\Bbb R }

\def\Z{\Bbb Z }
\def\C{\Bbb C }

%
%
\font\bigbf=cmbx10 scaled 1200
\font\author=cmcsc10 scaled 1200

\font\small=cmr8
\TagsOnRight
\NoRunningHeads
\def\np{\vfil\eject}
\def\nl{\hfil\newline}

\redefine\l{\lambda}

\def\w{\omega}
\def\a{\alpha}
\redefine\b{\beta}
\redefine\t{\tau}
\redefine\i{{\,\bold i\,}}
\def\PL{Phys\. Lett\. B}

\def\LMP{Lett\. Math\. Phys\. }

\def\FA{Funktional Anal. i. Prilozhen.}

\def\cintu{\frac 1{2\pi\i}\oint_{C_{u}}}

\def\im{\text{Im\kern1.0pt }}
\def\re{\text{Re\kern1.0pt }}
\def\res{\text{res}}
\def\RS{Riemann surface}

\def\ord{\operatorname{ord}}
\def\fpz{\frac {d }{dz}}

\def\pfz#1{\frac {d#1}{dz}}

\def\KNN {Krichever - Novikov }

%
%
\hfill LMU-TPW 92-05

\hfill Mannheimer Manuskripte 141

\hfill May 1992

\vskip 1cm
\centerline{\bigbf String Branchings on Complex Tori and Algebraic }
\vskip 0.3cm
\centerline{\bigbf Representations of Generalized \KNN Algebras}
\vskip 1.5cm
\centerline{\author Andreas Ruffing${}^{1,2},$\qquad
\qquad Thomas Deck ${}^{3}$,}
\vskip 0.3cm
\centerline{\author
Martin Schlichenmaier${}^{3}$
}
\vfill
\def\adressem{
{\baselineskip=10pt plus 0.2pt
${}^1$\small Sektion Physik der Universit\"at M\"unchen,
Theresienstra\ss e 37,
D-W-8000 M\"unchen 2, Germany}
}
\def\adressek{
{\baselineskip=10pt plus 0.2pt
${}^2$\small Institut f\"ur Theoretische Physik,
Universit\"at Karlsruhe,
D-W-7500 Karlsruhe 1, Germany
}}
\def\adressema{
{\baselineskip=10pt plus 0.2pt
${}^3$\small Fakult\"at f\"ur Mathematik und Informatik,
Universit\"at Mannheim, A5, Postfach 103462,   \nl
D-W-6800 Mannheim 1, Germany,
}}
\noindent\adressem \vskip 0.2cm
\noindent\adressek \vskip 0.2cm
\noindent\adressema
\vskip 1cm\noindent
{\bf Abstract.}
The propagation
differential for bosonic strings on a complex torus
with three symmetric punctures is investigated. We study
deformation aspects between two point  and three point differentials
as well  as the behaviour of the corresponding
Krichever-Novikov algebras.
The structure constants are calculated and from this we derive a
central extension of the Krichever-Novikov algebras by means
of $b-c$ systems.
The defining cocycle for this central extension deforms to
the well known Virasoro cocycle for certain kinds of
degenerations of the torus.
\np\pageno=1
{\bf 1. Introduction}
\vskip 0.5cm
String theory and conformal field theories were starting points for the
investigation of Lie algebras of meromorphic vector fields on
Riemann surfaces in the context of mathematical physics. In 1987
Krichever and Novikov \cite{11} studied (beside other objects)
algebras of meromorphic vector fields with poles only at
fixed points $P_+,P_-$ (called markings or punctures)
 on a compact \RS\  $X$,
as well as central extensions of such algebras, see also \cite{14}.
These algebras are nowadays called Krichever-Novikov algebras
(KN algebras).

The generalization of their concept for more than two markings of a
surface was given by Dick \cite{4-8} and Schlichenmaier \cite{15-18}.
They studied algebras of meromorphic vector fields on a genus $g$
\RS\  which have poles only at a finite number of fixed points.
For related work see also \cite{9}.
In \cite{2,3} it was shown by Deck that a special choice of
two punctures on a torus provides simple structure constants.
They only contain the parameters of the Weierstra\ss\ $\wp$--function
$$
e_1:=\wp(\frac 12),\qquad
e_2:=\wp(\frac 12+\frac {\tau}2),\qquad
e_3:=\wp(\frac {\tau}2)\ .\tag 1$$
Here $\tau$ denotes the normalized period of the torus.
This fact allows one to consider deformations of the analytic
structure in terms of $\ e_1,e_2,e_3$. In order to obtain similar
results on a complex torus with three punctures (markings) and
in order to study some further aspects of deformation
(concerning the number of punctures) let us repeat some
basic definitions.

Let $X$ be a compact \RS\ of genus $g$ with $N$ fixed points ($N\ge 2$).
We denote the set of these points with $S$. Let $S$ split into two nonempty
disjoint subsets $E$ and $O$, $S=E\cup O$. $E$ is called the set of
``in-points'' and $O$ the set of ``out-points''.
Let $\  \omega\ $ be the unique meromorphic differential with only
imaginary periods, holomorphic on $X\setminus S$ with poles of order
1 at $S$ and with the residues
$$\res_P(\omega)=+\frac 1{\#E},\quad P\in E,\qquad\text{and}\qquad
\res_Q(\omega)=-\frac 1{\#O},\quad Q\in O\ .\tag 2$$
In the following we restrict our investigations to $g=1$
(i.e. to the case of complex tori) and the
special cases of one  in-point $\ P_1\ $ and two out-points $\ Q_1,Q_2\ $.
The well known case of two points will appear as a degenerated
three point case.
The generalizations to more than 3 points are discussed in \cite{13}.

We  choose a point $R$ on $X\setminus S$.
Due to the above properties the function
$$t(P):=\re\int_R^P\omega\tag 3$$
is  well defined and harmonic on $X\setminus S$.
The level lines
$\ C_u:=\{P\in X\mid t(P)=u\}\ $ for $u\in\R$
define a global fibration of the surface $X\setminus S$ interpreted
as the ``position'' of the string(s) at ``time'' $u$.
For large negative $u$
these lines are small circles around the points in $E$
and for large positive $u$ they are small circles around the points in
$O$. So the strings ``travel'' from the points with positive
residues to those with negative residues.
They are called ``ingoing'' strings at $\ E\ $ and ``outgoing'' strings
at $\ O\ $.
At the points were $\omega$ has a zero the level lines (strings)
split or rejoin. Because of this physical interpretation
$\omega $ is called propagation differential or
differential of the time development. A physical realization of the
KN algebra (with central extensions) is
given by the commutation relations of the
components of the energy-momentum tensor which belong to a string
moving in Minkowski space \cite{11}.
\vskip 0.5cm
{\bf 2. The propagation differential}
\vskip 0.5cm
We consider a torus $T=\C/ L$ with
$\  L:=\{\ z\in\C\mid z=m+n\t,\quad m,n\in\Z\ \}$
\nl and $\ \t\in\C,\ \im \t>0\ $.

We need the following facts about the Weierstra\ss\ $\wp$--function
\cite{10}:\nl
$\wp$ is an even meromorphic doubly periodic function (i.e.
$\wp(z+w)=\wp(z),\ \forall w\in  L$)
with a second order pole only in the lattice points and
holomorphic elsewhere with Laurent expansion
$\ \wp(z)=z^{-2}+\cdots.$ By the covering map
$\ \pi:\C\to\C/ L$ we have a holomorphic 1-form $dz$ on
the torus $T$.

As markings we choose
the following points
$$P_1=0 \mod  L,\qquad
Q_1=1/2+q \mod  L,\qquad
Q_2=1/2-q \mod  L,\tag 4$$
with $\ q\ne 0,\ 1/2 \mod L$.

\noindent We construct now a differential $\ \w\ $ with the properties

(1) $\w$ is holomorphic on $T\setminus\{P_1,Q_1,Q_2\}$,

(2) $\res_{P_1}(\w)=+1$\quad and\quad
$\res_{Q_i}(\w)=-1/2,\  i=1,2$,

(3) $\w$ has only imaginary periods.

\noindent From considerations about the (probable) zeros of $\w$
(the set $\{1/2,\t/2,(1+\tau)/2\}$) we found the following result
\proclaim{Proposition}
The propagation differential with respect to the punctures (4)
is given by
$$\w_q(z)=\widehat{\w}_q(z)\,dz\
 =\ -\frac 12\,\frac {\wp'(z)}{\wp(z)-\wp(1/2+q)}\,dz\ .\tag 5$$
\endproclaim
\demo{Proof}
The singularities of $\w$ (all of first order) follow from the facts
about $\wp(z)$ as well as the property
$\res_0(\w)=1$.
\nl Also we have $\ \wp(z)=-\wp(-z)\ $ and this
yields together with $\ \sum_{Q\in T}\res_Q(\w)=0\ $ (residue theorem)
$\ \res_{1/2+q}(\w)=\res_{1/2-q}(\w)=-1/2\ $
and hence the property (2).\nl
The final requirement (3) follows from the fact that
$\w$ is antisymmetric in the argument with respect to
the  points $\ 0,\ 1/2,\ (1+\tau)/2\ $ and $\ \tau/2\ $
and has purely real residues.\qed
\enddemo
\noindent Note, in the proof above we used $\ \res_{z \mod L}=\res_z\ $.
\vskip 0.5cm
{\bf 3. Degenerated case and calculation of the separation time}
\vskip 0.5cm
If we now set $q=0$ we will find a holomorphic differential on a
2-punctured torus:
$$\w_0(z)=\widehat{\w}_0(z)\,dz
=-\frac 12\,\frac {\wp'(z)}{\wp(z)-e_1}\,dz\ .\tag 6$$
Obviously we have $\res_0(\w_0)=-\res_{1/2}(\w_0)=+1$. All other required
properties of a propagation differential are conserved and so $\w_0$
describes a process of pure generation and annihilation of a
 string-pairing.
In the limit $q=0$ the two poles
 of $\ \w_q\ $ for outgoing strings coincide in $1/2$. This
case has been studied in \cite{3} and is a special case in our
prescription. The differential $\w_q$ has a big advantage:
It's real part can be integrated at once because it has
the structure of a logarithmic differential.
If we integrate between the two zeros
$\ \tau/2\ $ and $\ (1+\tau)/2\ $ of $\ \w_q\ $ we find
$$(\Delta t)_q:=\re\int_{\t/2}^{(1+\t)/2}\w_q=\frac 12
\ln\left|\frac {e_3-\wp(1/2+q)}{e_2-\wp(1/2+q)}\right|$$
and
$$(\Delta t)_0=\frac 12
\ln\left|\frac {e_3-e_1}{e_2-e_1}\right|\tag 7$$
as an obvious limit for $q\to 0$.

{}From the physical background we call $\ (\Delta t)_q\ $ separation time.
We only give an interpretation of (7). Here we find
$\ (\Delta t)_0=-\frac 12\ln|\mu|\ $ where the parameter
\nl $\mu:=(e_2-e_1)(e_3-e_1)^{-1}\ $
is directly related to the algebraic equation of $\wp,\wp'$.

More precisely, $\mu\ $ is the moduli parameter for tori with
level 2
structure, i.e. for tori with fixed 2-torsion points.
If we introduce the principal congruence subgroup of level 2
$$\Gamma(2):=\{A\in \text{SL}(2,\Z)\mid
\ A\equiv \pmatrix 1&0\\ 0&1\endpmatrix\mod 2\ \}$$
of the full modular group $\ \text{SL}(2,\Z)\ $ then the moduli space
of tori with level 2 structure is the quotient of the
upper half-space $\ \Cal H:=\{z\in\C\mid \im z>0\}\ $ with respect to
the action of $\ \Gamma(2)\ $ by fractional linear transformation.
In other words
$$z\sim z'\iff z'=\frac {az+b}{cz+d},\qquad
\pmatrix a&b\\c&d\endpmatrix\in\Gamma(2)\ .\tag 8 $$
It is
well known (see \cite{12, App.7}) that $\ |\mu|=1\ $ iff
the lattice parameter $\ \tau\ $ is equivalent to a $\ \tau'\ $ under
the relation (8) with $\ \re \tau'=\pm 1/2\ $.
\vskip 0.5cm
{\bf 4. The corresponding Lie  algebra of vector fields}
\vskip 0.5cm
The KN algebra is given by those meromorphic vector fields which are
holomorphic on $T\setminus S$. There is a natural identification
of meromorphic functions and meromorphic vector fields on a
complex torus
by the correspondence $\ f(z)\mapsto f(z)\fpz\ \ $.
Let us therefore consider an infinite
set $\ \Cal A=\{A_k\}_{k\in\Z}\ $ of functions.
The elements are
given by
$$\align
A_n(z)&:=\exp(n\int_{\l}^z\w_q)=a_n(\l)\left(\wp(z)-\wp(1/2+q)
\right)^{-n/2},\quad n\in 2\Z,\tag 9\\
A_{\a}(z)&:=\widehat{\w}_q(z)A_{\a+1}(z),\qquad\qquad\a\in 2\Z+1 \ .\tag 10
\endalign
$$
We choose $\ \l\in\C\ $ such that $\ a_2(\l)=1\ $ and hence all
$a_n(\l)=1$. Note that $A_n=(A_2)^{n/2}$. It is worth-while
mentioning that all these functions are given in terms of
the propagation differential $\w_q$.
Obviously the $A_n$ are even functions whereas the $A_{\a}$ are odd
functions. From now on we denote even indices with $m,n$ and odd ones
with $\a,\b$. For convenience let us introduce an order triple
of a function $f$:
$$f^*:=\ (\ord_0f,\ord_{1/2+q}f,\ord_{1/2-q}f)\ \in\ \Z^3\ .\tag 11$$
We then have
$$A_n^*=(n,-n/2,-n/2)\ \qquad\text{and}\qquad
 A_{\a}^*=(\a,(-\a-3)/2,(-\a-3)/2)\ .\tag 12$$
Considering the $\ \ord_0\ $ entry  it is  evident that
the set $\ \Cal A\ $ is linearly independent over $\C$.
\proclaim{Proposition}
A basis of meromorphic functions on $T$ which are holomorphic on
$T\setminus S$ is given by the set $\Cal A$.
\endproclaim
\demo{Proof}
There are at least two ways to show that $\Cal A$ represents a basis
of vector fields. Each way gives insights which will we useful
later on.

1. One can show that there exists linear combinations
$E_1,E_2,E_3,E_4$ by elements of $\Cal A$ so that
$\ E_1^*=(-1,-1,0)\ $,
$\ E_2^*=(0,-1,-1)\ $,
$\ E_3^*=(0,0,2)\ $,
$\ E_4^*=(2,-1,-1)\ $.
For instance
such an element $E_1$ can be obtained
as a certain
linear combination of $A_2$ and $A_{-1}$
by looking at their analytic behaviour
near $\ 1/2\pm q\ $.
It can be shown that any function $f$ on the torus with markings in $S$
can be algebraically combined by $E_1,E_2,E_3,E_4$ and $A_0$. We have
$$ A_iA_j=A_jA_i=A_{i+j}\quad\text{for}\quad i\in 2\Z,\ j\in\Z\ \tag 13$$
and $\ (\widehat{\w}_q^2)^*=(-2,-2,-2)$. As
$\ \widehat{\w}_q,\ A_{-2},\ A_0,\ A_2,\ A_4\ $
are even function one can show using (12)
that
$$\widehat{\w}_q^2=\l_1A_{-2}+\l_2A_{0}+\l_3A_{2}+\l_4A_{4},
\quad\text{with}\quad
\l_1,\l_2,\l_3,\l_4\in\C\ .\tag 14$$
Hence, we obtain with the same values for $\l_i$
$$ A_iA_j=
\l_1A_{i+j}+\l_2A_{i+j+2}+\l_3A_{i+j+4}+\l_4A_{i+j+6}
\quad\text{for}\quad i,j\in 2\Z+1\ .\tag 15$$
Now one can easily see that any $f$ is a linear combination of
elements of $A$.

2. Because $A_n$
is even and $A_{\a}$ is odd there are linear combinations
\nl $\ B_k=\gamma A_{2k}+\delta A_{2k-3},\ $
$\ C_k=\gamma' A_{2k}+\delta' A_{2k-3}\  $
with order triple
$\ B_k^*=(2k-3,-k,b)$,\nl $C_k^*=(2k-3,c,-k)$ and $b,c>-k$.
By a simple elimination argument concerning the pole orders of
an arbitrary function $f$ ( of the prescribed type) it is
verified that $f$ is a linear combination of the
$A_k,B_k,C_k$ and consequently a linear combination of the set $\Cal A$
alone.\qed
\enddemo
The next step is  the computation of the Lie algebra of vector fields
which are given by $\ l_k(z):=A_k(z)\pfz\ $ under Lie bracket.
{}From (9) and (10) we find
$$\align
[l_n,l_m](z)&=(m-n)\;l_{m+n-1}(z),\tag 16\\
[l_\a,l_\b](z)&=(\b-\a)\;\widehat{\w}_q^2(z)\;l_{\a+\b+1}(z),\tag 17\\
[l_\a,l_n](z)&=\left((n-\a+1)\;\widehat{\w}_q^2(z)+
\widehat{\w}'_q(z)\right)l_{n+\a+1}(z)\ .
\tag 18
\endalign$$
It is an interesting fact that the whole structure of the algebra
is encoded in $\w_q$ \cite{13}.
The element $\widehat{\w}_q'$ is
(like $\widehat{\w}_q^2$) an even functions hence a
linear combination of the even elements $A_n$. To obtain its explicit
representation we use $\ A_n'(z)=n\widehat{\w}_q(z)A_n(z)\ $. Using this and
the representation (14) for $\widehat{\w}_q^2$ we can deduce from
$\ (\widehat{\w}_q^2)'(z)=2\widehat{\w}_q(z)\widehat{\w}_q'(z)\ $
$$\widehat{\w}_q'(z)=-2\l_1A_{-2}(z)+ 2\l_3A_{2}(z)
+4\l_4A_{4}(z) \ .\tag 19$$
Together with property (13) the structure coefficients of the
algebra immediately follow from (16)--(18). The remaining task is to
find the coefficients $\l_1,\l_2,\l_3,\l_4$.
But this can be done in an algebraic way:

Substituting the Weierstra\ss\ relation
$\ \wp'=4(\wp-e_1)(\wp-e_2)(\wp-e_3)\ $ into
$$
\widehat{\w}_q^2(z)=\frac 14\, \frac {\wp'{}^2(z)}
{\left(\wp(z)-\wp(1/2+q)\right)^2}
\tag 20$$
as well as the expression (9) into (14) we find a polynomial identity in the
variable $X:=\wp(z)$ from which we can recognize the coefficients
$\l_1,\l_2,\l_3,\l_4$.
The computation finally yields the required generalized
Krichever-Novikov algebra:
$$\align
[l_n,l_m]\ &=\ (m-n)\;l_{n+m-1},\tag 21\\
[l_\a,l_\b]\ &=\ (\b-\a)\big(l_{\a+\b-1}+3\,\wp(1/2+q)\;l_{\a+\b+1}+\\
&\qquad (3\wp^2(1/2+q)-(e_2^2+e_2e_3+e_3^2))\;l_{\a+\b+3}+\\
&\qquad 1/4\,\wp'{}^2(1/2+q)\;l_{\a+\b+5}\big)\ ,\tag 22\\
[l_\a,l_n]\ &=\ (n-\a)l_{\a+n-1}+(n-\a-1)3\wp(1/2+q)\;l_{\a+\b+1}+\\
&\qquad (n-\a-2)(3\wp^2(1/2+q)-(e_2^2+e_2e_3+e_3^2))\;l_{\a+n+3}+\\
&\qquad (n-\a-3)\, 1/4\,\wp'{}^2(1/2+q)\;l_{\a+n+5}\ .\tag 23
\endalign
$$
Indeed the case $q=0$ coincides with that given in \cite{3}.
In this sense we have another aspect of deformation in the
family of KN algebras. The two point algebra ($q=0$) is
obtained by taking the limit $\ q\to 0\ $
of the three point algebra:
$$\align
[l_n,l_m]\ &=\ (m-n)\;l_{n+m-1},\tag 24\\
[l_\a,l_\b]\ &=\ (\b-\a)\big(l_{\a+\b-1}+3e_1l_{\a+\b+1}+
(e_1-e_2)(e_1-e_3)\;l_{\a+\b+3})
\ ,\tag 25\\
[l_\a,l_n]\ &=\ (n-\a)\;l_{\a+n-1}+(n-\a-1)3e_1l_{\a+\b+1}+\\
&\qquad (n-\a-2)(e_1-e_2)(e_1-e_3)\;l_{\a+n+3}
\ .\tag 26
\endalign
   $$
The algebraic-geometric background for further degenerations has
been studied in detail by one of us \cite{19}. Note for
instance that by $\ e_1=e_2=e_3=0\ $ this algebra becomes a
representation for the centerless Virasoro algebra.
\loadbold
\vskip 0.5cm
{\bf 5. Central extensions by $\boldkey b \boldkey -\boldkey c$ systems}
\vskip 0.5cm
Such central extensions in general for two point cases were
studied in \cite{1}.
The properties of $\ b-c\  $systems have been generalized to more than two
punctures
\cite{17,18}. Explicit expressions for our case of two points had been
obtained by the formula given in \cite{1}.
In this last section we investigate
some aspects of these generalizations to the algebra (21)--(23).
In order to get expressions which are easier to compare with the centrally
extended Virasoro algebra we perform an index shift
$\ e_i:=
l_{i+1}$ and denote by $C_{ij}^k$ the structure constants of
$\ [e_i,e_j]$.

By integrating the product of a $\l-$form and a $(1-\l)-$form over a
nonsingular level line $C_u$ a dual pairing is introduced.
It does not depend on $u\in\R  $ (see \cite{17} for details).
Here we consider only $(\l,1-\l)=(-1,2)$.
\proclaim{Proposition}
The objects $\ \Omega_k:=A_{-k-2}(dz)^2\ $ and $e_j$ with
$k,j\in\Z$ form a dual system by the product
$$\langle e,\Omega\rangle:=\cintu e\cdot\Omega\ .\tag 27$$
\endproclaim
\demo{Proof}
Because all curves $C_u$ are homologous one can use
local residue calculus around $P_1$ as well as around $Q_1$ and $Q_2$.
These two integrations (residue calculus) give two conditions
for the non-vanishing of
$\ \langle e_j,\Omega^k\rangle\  $
which taken together yield
$$ \text{i}(e_j)(\Omega_k):=\ \langle e_j,\Omega^k\rangle=c_k\delta_i^k\ .$$
$c_k=1$ follows from the exact calculation of the residue at $P_1$.
\qed\enddemo
Let $V$  be the vector space generated by all semi-infinite forms
of weight 2
$$\Omega^{i_1}\wedge
\Omega^{i_2}\wedge\cdots\wedge\Omega^{s}\wedge
\Omega^{s-1}\wedge\Omega^{s-2}\wedge\cdots,\tag 28$$
with $\ i_1>i_2>\ldots s>s-1\ldots\ $ and
where the sequence of $\Omega^i$ contains all indices smaller
than $s$ (for some $s$ which is not fixed).
Furthermore we have operators
$\ c^i,\ b_k\ $ defined by
$\ c^i:=\Omega^i\wedge\ $ and
$\ b_k:=\text{i}(e_k)$ + product rule (anticommuting).
As one can verify these operators satisfy a Clifford algebra
$$\{b_k,c^i\}:=b_k\circ c^i+c^i\circ b_k=\delta_k^i,\qquad
\{b_k,b_l\}=\{c^i,c^j\}=0\ .\tag 29$$
The following results are a concrete application
of \cite{1} and \cite{18}.
The vacuum vector
$\ |\;0\,\rangle:=\Omega^{-2}\wedge\Omega^{-3}\wedge\Omega^{-4}\wedge
\cdots\ $ has the properties
$$ c^i|\;0\,\rangle =0,\quad i<-1\qquad \text{and}\quad
b_k|\;0\,\rangle =0,\quad k\ge -1\ .$$
We get a normal ordering by $\ :b_kc^i:\ := b_kc^i\ $ for $i<-1$
and $\ -c^ib_k\ $ for $\ i\ge -1$.

For  general  $b-c$ systems of weight $(\lambda,1-\l)$ the
energy momentum tensor  is defined by \cite{18}
$$T(z):=\ :(1-\l)c(z)\partial_zb(z)-\l (\partial_z c(z))b(z):\ .\tag 30$$
$T$ is a form of weight $2$ and can be expanded as
$\ T(z)=\sum_kL_k\Omega^k(z)\ $.
The operator valued coefficents are
$$L_i=\cintu T\cdot e_i$$
by the duality (27).
Here we have
$\l=2$. Using the expansion
$c(z)=\sum_ic^ie_i(z)$, $b(z)=\sum_j b_j\Omega^j(z)$
and recognizing (30) essentially as the Lie derivative
\cite{18}
$$L_c(g_\l)=c\partial_zg_\l+\l g_\l\partial_zc\
\quad (g_\l \text{ a $\l$-form})$$
we finally obtain the result
$$L_i=\sum_{j,k}C_{ij}^k :b_kc^j:\ .$$
By explicit calculation one can see  that these elements satisfy a
centrally extended KN algebra for the singularity set
$S$ and cause the relation
$$[L_i,L_j]=-\sum_{k}C_{ij}^kL_k+\chi_{ij}\quad
\text{with}\quad
\chi_{ij}=\left(\sum_A-\sum_B\right)C_{ik}^lC_{jl}^k\ .\tag 31$$
The sets $A$ and $B$ of  double indices are given by
$\ A:=\{k<-1,\ l\ge -1\}\ $,\nl
$\ B:=\{k\ge -1,\ l< -1\}\ $.
It is worthwhile noting that this formula (31) is independent
from the number of markings.
The only fact which is important with respect to this kind of
representations is the generalized grading of the algebra. This
grading is responsible for the $\ L_i\ $ to be well defined
operators on $V$  because
then $\ L_i\cdot $(basis vector)\  always is a finite sum.

Let us give the cocycle and show its formal deformation to the cocycle
of the Virasoro algebra. By direct but tedious calculations
one obtains the following results (see \cite{13} for details):
First we get
$$
\chi_{\a n}=\chi_{n\a}=0\quad\text{for}\quad
n,\a+1\in 2\Z\ .\tag 32$$
To write down the remaining terms we use the following
abbreviations:
$\ I:=\{4,5,6,7\}\ $, $\b_i:=\b+i$, $m_i=m+i$  and
$$\gathered \l_4:=1,\quad \l_5:=3\;\wp(1/2+q),\\
\l_6:=3\wp^2(1/2+q)-(e_2^2+e_2e_3+e_3^2),\quad
\l_7:=1/4\wp'{}^2(1/2+q)\endgathered\tag 33$$
for the constants coming from the structure constants $C_{kl}^r$.
In addition we define
$$Q_k:=\sum_{i,j\in I}\l_i\l_j\delta_{i\cdot  j}^k\ .\tag 34$$
With this we receive (see \cite{13})
$$\align
\chi_{\a\b}\quad&=\
13/6\;(\b^3-\b)\delta_{\a+\b}^0+
(13/6\;\b_1^3-2/3\;\b_1)\delta_{\a+\b}^{-2}Q_{20}+
\\&\qquad
(13/6\;\b_2^3-25/6\;\b_2)\delta_{\a+\b}^{-4}Q_{24}+
(13/6\;\b_3^3-76/6\;\b_3)\delta_{\a+\b}^{-6}Q_{28}^*,\tag 35
\endalign$$
$$\align
\chi_{nm}\quad&=\
13/6\;(m^3-m)\delta_{m+n}^0+
(13/3\;m_1^3+5/3\;m_1)\delta_{m+n}^{-2}Q_{20}+
\\&\qquad
(13/3\;m_2^3+11/3\;m_2)\delta_{m+n}^{-4}Q_{24}+
(13/6\;m_2^3-2/3\;m_2)\delta_{m+n}^{-4}Q_{25}+
\\&\qquad
(13/3\;m_3^3+5/3\;m_3)\delta_{m+n}^{-6}Q_{28}^*+
(13/3\;m_3^3-25/3\;m_3)\delta_{m+n}^{-6}Q_{30}+
\\&\qquad
(13/3\;m_4^3-58/3\;m_4)\delta_{m+n}^{-8}Q_{35}^*+
(13/6\;m_4^3-73/6\;m_4)\delta_{m+n}^{-8}Q_{36}+
\\&\qquad
(13/3\;m_5^3-133/3\;m_5)\delta_{m+n}^{-10}Q_{42}^*+
(13/6\;m_6^3-110/3\;m_6)\delta_{m+n}^{-12}Q_{49}^*\ .
\tag 36
\endalign
$$
Here we have denoted all those $Q_k$ with a star which
will  vanish for $q=0$.

In a first step we degenerate to the two point case.
Because the structure constants depend continuously
 on $q$ and because they
reduce to those in the two-point case if $q$ equals zero we necessarily
get the cocycles of the two-point case.
In a second step we degenerate
to $\ e_1=e_2=e_3=0\ $ and find that all remaining terms in
which $Q_k$ occur become zero. What remains ist the well known
extension form the Witt algebra to the Virasoro algebra
$$[L_m,L_n]=(m-n)L_{m+n}+\frac {13}6(m^3-m)\delta_{m+n}^0,\quad m,n\in\Z\ .$$
This expected result is a further motivation for the research
on degenerated tori.
\vskip 0.3cm
{\bf Acknowledgements}
\vskip 0.3cm
\noindent The first author (A.R.) likes to thank
R.~Dick for very fruitful discussions.

\bye